\documentclass
[12pt]
{amsart}

\usepackage{latexsym,amsmath}

\theoremstyle{plain}
  \newtheorem{thm}{Theorem}[section]
  \newtheorem{prop}[thm]{Proposition}
  \newtheorem{lem}[thm]{Lemma}
  \newtheorem{cor}[thm]{Corollary}
  \newtheorem{conj}[thm]{Conjecture}
\theoremstyle{definition}
  \newtheorem{dfn}[thm]{Definition}
  \newtheorem{exmp}[thm]{Example}
  \newtheorem{prob}[thm]{Problem}
\theoremstyle{remark}
  \newtheorem{rem}[thm]{Remark}

\newcommand{\Ass}{\operatorname{Ass}}
\newcommand{\codim}{\operatorname{codim}}
\newcommand{\depth}{\operatorname{depth}}

\newcommand{\hull}{\operatorname{hull}}
\newcommand{\init}{\operatorname{in}}
\newcommand{\interior}{\operatorname{int}}
\newcommand{\irr}{\operatorname{Irr}}
\newcommand{\lcm}{\operatorname{lcm}}
\newcommand{\pdim}{\operatorname{proj-dim}_S}
\newcommand{\rad}{\operatorname{rad}}
\newcommand{\rank}{\operatorname{rank}}
\newcommand{\Tor}{\operatorname{Tor}}

\newcommand{\supp}{\operatorname{supp}}

\def\aa{{{\rm \bf a}}}
\def\bb{{{\rm \bf b}}}
\def\cc{{{\rm \bf c}}}

\def\mm{{\mathfrak m}}
\def\xx{{\bf x}}
\def\DD{{\bf D}}
\def\LL{{\mathcal L}}
\def\NN{{\mathbb N}}

\def\ZZ{{\mathbb Z}}
\def\mmm#1#2{m_{\{#1,#2\}}}
\def\<{{\langle}}
\def\>{{\rangle}}

\begin{document}
\title{Generic and cogeneric monomial ideals}
\author{Ezra Miller}
\address{Department of Mathematics, 
University of California, 
Berkeley, CA  94720 USA}
\email{enmiller@math.berkeley.edu}

\author{Bernd Sturmfels} 
\address{Department of Mathematics, 
University of California, 
Berkeley, CA 94720 USA}
\email{bernd@math.berkeley.edu}
                                           \
\author{Kohji Yanagawa}
\address{Graduate School of Science, Osaka University, Toyonaka, Osaka 
560, Japan}
\email{yanagawa@math.sci.osaka-u.ac.jp}

\begin{abstract}
Monomial ideals which are generic with respect to either their generators
or irreducible components have minimal free resolutions derived from
simplicial complexes.  For a generic monomial ideal, the associated
primes satisfy a saturated chain condition, and the Cohen-Macaulay
property implies shellability for both the Scarf complex and the
Stanley-Reisner complex.  Reverse lexicographic initial ideals of generic
lattice ideals are generic.  Cohen-Macaulayness for cogeneric ideals is
characterized combinatorially; in the cogeneric case the Cohen-Macaulay
type is greater than or equal to the number of irreducible components.
Methods of proof include Alexander duality and Stanley's theory of local
$h$-vectors.
\end{abstract}

\maketitle

\section{Genericity of Monomial Ideals Revisited}
\label{generic}

Let $M$ be a monomial ideal minimally generated by monomials $m_1,
\ldots, m_r$ in a polynomial ring $S = k[x_1, \ldots, x_n]$ over a field
$k$.  For a subset $\sigma \subseteq \{1, \ldots, r \}$, we set
$m_{\sigma} := \lcm ( m_i \, | \, i \in \sigma )$, and
$\aa_\sigma := \deg m_{\sigma} \in {\mathbb N}^n$ the exponent vector of
$m_{\sigma}$.  Here $m_{\emptyset} = 1$.  For a monomial $\xx^\aa =
x_1^{a_1} \cdots x_n^{a_n}$, we set $\deg_{x_i} (\xx^\aa) := a_i$, and we
call $\supp (\xx^\aa) := \{ i \, | \, a_i \ne 0\} \subseteq \{1, \ldots,
n\}$ the {\it support} of $\xx^\aa$.

\begin{dfn} \label{gendef} A monomial ideal
$M= \<m_1, \ldots, m_r\>$ is called {\it generic} 
if for any two distinct generators $m_i, m_j $ of $M$
which have the same positive degree in some
variable $x_s$ there 
exists a third monomial generator $m_l \in M$ which divides $m_{\{i,j\}} = \lcm
(m_i, m_j)$ and satisfies $\supp(m_{\{i,j\}}/m_l) = \supp(m_{\{i,j\}})$.
\end{dfn}

The above definition of genericity is more inclusive than the one given
by Bayer-Peeva-Sturmfels \cite{BPS}, but we will see that this definition
permits the same algebraic conclusions as the one in \cite{BPS}.  There
are important families of monomial ideals which are generic in the sense
of Definition~\ref{gendef} but not in the sense of \cite{BPS}.  One such
family is the initial ideals of generic lattice ideals as in Theorem
\ref{ingen}.  Here is another one:

\begin{exmp} \label{tree}
The {\it tree ideal} $\,M = \<\,\bigl(\prod_{s \in I}
x_s\bigr)^{n-|I|+1} \mid \emptyset \neq I \subseteq \{1,\ldots,n\}\>\,$
is generic in the new sense but very far from generic in the old sense. 
This ideal is Artinian of colength $(n+1)^{n-1}$, the 
number of trees on $n+1$ labelled vertices. 
\end{exmp}

Recall that a monomial ideal $M \subset S$ can be uniquely written as a
finite  irredundant intersection $M = \bigcap_{i=1}^r M_i$ 
of irreducible monomial ideals 
(i.e., ideals generated by powers of variables).  
We say $M_i$ is an {\it irreducible component} of $M$. 

\begin{dfn}\label{cogendef}
A monomial ideal with irreducible decomposition $\, M =
\bigcap_{i=1}^r M_i\,$  is  called {\it cogeneric} if the following condition
holds: if distinct irreducible components $M_i$ and $M_j$ have a
minimal generator in common, there is an irreducible component $M_l
\subset M_i + M_j$ such that $M_l$ and $M_i + M_j$ do not have a minimal
generator in common. 
\end{dfn}

A monomial ideal $M$ is cogeneric if and only if its {\it Alexander dual} 
$M^\aa$ is generic. See \cite{Mil} or Section 4 for the
relevant definitions.
Cogeneric monomial ideals will be studied in
detail in Section 4. 
The remainder of this section is devoted to
basic properties of generic monomial ideals.   

Let $M \subset S$ be a monomial ideal minimally generated by monomials
$m_1, \ldots, m_r$ again.  The following simplicial complex on $r$
vertices, called the {\it Scarf complex} of $M$, was introduced by Bayer,
Peeva and Sturmfels in \cite{BPS}:
$$
  \Delta_{M} \quad :=\quad \{ \sigma \subseteq \{1, \ldots, r \} \, | \, 
  \text{$m_{\tau} \ne m_\sigma$ for all $\tau \ne \sigma$} \}.
$$

Let $S(-\aa_\sigma)$ denote the free $S$-module with one generator
$e_{\sigma}$ in multidegree $\aa_{\sigma}$.  The {\it algebraic Scarf
complex} $F_{\Delta_M}$ is the free $S$-module $\bigoplus _{\sigma \in
\Delta_M} S(-\aa_\sigma)$ with the differential
$$  d(e_{\sigma}) \,\,\, =  \,\,\,
\sum_{i \in {\sigma}} \operatorname{sign}(i,\sigma)
  \cdot \frac{m_{\sigma}}{m_{\sigma \setminus \{i\}}} \cdot e_{\sigma
  \setminus \{i\}}
$$
where $\operatorname{sign}(i,\sigma)$ is $(-1)^{j+1}$ if $i$ is the
$j$-th element in the ordering of $\sigma$.  It is known that
$F_{\Delta_M}$ is always contained in the minimal free resolution of
$S/M$ as a subcomplex \cite[\S 3]{BPS}, although $F_{\Delta_M}$ need not
be acyclic in general.  However we will see in Theorem~\ref{new-def} that
it is acyclic if $M$ is generic, as was the case under the old
definition.

\begin{lem}\label{1st_lemma}
Let $M = \<m_1, \ldots, m_r\>$ be a generic monomial ideal.  If $\sigma
\not \in \Delta_M$, then there is a monomial $m \in M$ such that $m $
divides $m_{\sigma}$ and $\supp (m_{\sigma} /m) = \supp (m_{\sigma})$.
\end{lem}

\begin{proof}
Choose $\sigma \not\in \Delta_M$ maximal among subsets of $\{1,\ldots,r\}$
with label $\aa_\sigma$. Then
$m_{\sigma} = m_{\sigma \setminus \{ i\}}$ for some $i \in \sigma$.  If
$\supp(m_{\sigma}/m_i) = \supp (m_{\sigma})$, the proof is done.
Otherwise, there is $\sigma \ni j \ne i$ with $\deg_{x_s} m_i =
\deg_{x_s} m_j > 0$ for some $x_s$.  Since $M$ is generic, there is a
monomial $m \in M$ which divides $m_{\{i,j\}}$ and satisfies $\supp(m_{\{
i,j\}}/m) = \supp (m_{\{i,j\}})$.  Since $m_{\{i,j\}}$ divides
$m_{\sigma}$, the monomial $m$ has the desired property.
\end{proof}

\vskip .5mm

The following theorem extends results in \cite{BPS} and is the main
result in this section.

\begin{thm}\label{new-def}
A monomial ideal $M$ is generic if and only if 
the following two  hold:

(a) The algebraic Scarf complex $F_{\Delta_M}$ equals the minimal free
	 resolution of $S/M$.

(b) No variable $x_s$ appears with the same non-zero exponent in $m_i$
	and $m_j$ for any edge $\{ i,j \}$ of the Scarf complex
	$\Delta_M$.
\end{thm}

\begin{proof}
Suppose that $M$ is generic. Then (b) is straightforward from the
definition, and, using Lemma~\ref{1st_lemma}, (a) is proved by the same
argument as in \cite[Theorem~3.2]{BPS}.

Assuming (a) and (b), we show that $M$ is generic.  For any generator
$m_i$ let 
$$  A_i \,\,:=\,\, \{m_j \mid m_j \neq m_i\ {\rm and}\ \deg_{x_s}
  m_j = \deg_{x_s} m_i > 0\ {\rm for\ some}\ s\}.
$$
The set $A_i$ can be partially ordered by letting $m_j \preceq m_{j'}$ if
$\mmm ij$ divides $\mmm i{j'}$.  It is enough to produce a monomial $m_l$
as in Definition~\ref{gendef} whenever $m_j \in A_i$ is a minimal element
for this partial order.  Supposing, then, that $m_j$ is minimal, use (a)
to write
\begin{equation} \label{syz}
  {\mmm ij \over m_i} \cdot e_i \,-\, {\mmm ij \over m_j} \cdot e_j
  \quad = \quad
  \sum_{\{u,v\} \in \Delta_M} b_{u,v} \cdot d(e_{\{u,v\}})
\end{equation}
where we may assume (by picking such an expression with a minimal number
of nonzero terms) that the monomials $b_{u,v}$ are $0$ unless $\mmm uv$
divides $\mmm ij$.  There is at least one monomial $m_l$ such that
$b_{l,j} \neq 0$, and we claim $m_l \not\in A_i$.  Indeed, $m_l$ divides
$\mmm ij$ because $\mmm lj$ does, so if $\deg_{x_t} m_i < \deg_{x_t} m_j$
(which must occur for some $t$ because $m_j$ does not divide $m_i$), then
$\deg_{x_t} m_l \leq \deg_{x_t} m_j$.  Applying (b) to $\mmm lj$ we get
$\deg_{x_t} m_l < \deg_{x_t} m_j$, and furthermore $\deg_{x_t} \mmm il <
\deg_{x_t} \mmm ij$, whence $m_l \not\in A_i$ by minimality of $m_j$.  So
if $\deg_{x_s} \mmm ij > 0$ for some $s$, then either $\deg_{x_s} m_l <
\deg_{x_s} m_j$ by (b), or $\deg_{x_s} m_l < \deg_{x_s} m_i$ because $m_l
\not\in A_i$.
\end{proof}

\begin{rem}
Condition (a) in Theorem \ref{new-def} splits into two parts: minimality
and acyclicity.  For the Scarf complex of {\it any} monomial ideal,
minimality is automatic since face labels $\aa_\sigma$ of $\Delta_M$ are
distinct.  It is acyclicity which must be checked.
\end{rem}

For an arbitrary monomial ideal $M$, Bayer and Sturmfels \cite[\S 2]{BS}
constructed a polyhedral complex $\hull(M)$ supporting a (not necessarily
minimal) free resolution of $M$.  Definition \ref{gendef} suffices to
imply that the hull complex equals the Scarf complex:

\begin{prop}\label{hull}
If $M$ is a generic monomial ideal, then the hull complex $\hull(M)$
coincides with $\Delta_M$, and in this case the hull resolution $F_{\hull
(M)} = F_{\Delta_M}$ is minimal.
\end{prop}

\begin{proof}
Essentially unchanged from the proof of \cite[Theorem~2.9]{BS}. 
\end{proof}

\noindent {\bf Example \ref{tree}} {\sl (continued) } The Scarf complex
$\Delta_M$ of $M$ is the first barycentric subdivision of the
$(n-1)$-simplex.  By Theorem~\ref{new-def}, $F_{\Delta_M}$ gives a
minimal free resolution of $S/M$.  Miller~\cite{Mil} also constructed a
minimal free resolution of $S/M$ as a {\it cohull resolution}, derived
essentially from the coboundary complex of a permutahedron.

\section{Associated Primes and Irreducible Components}
\label{assprime}

In this section we study the primary decomposition
 of a generic monomial ideal $M$.
For a monomial prime $P$ in $S$, we identify the homogeneous localization
$(S/M)_{(P)}$ with the algebra $k[x_i \mid x_i \in P] /M_{(P)}$, where
$M_{(P)}$ is the monomial ideal of $k[x_i \mid x_i \in P]$ gotten from 
$M$ by setting equal to $1\,$ all the variables not in $P$.

\begin{rem}
If $M$ is a generic monomial ideal then so is $M_{(P)}$.
\end{rem}

Let $M = \bigcap_{i=1}^r M_i$ be the irreducible decomposition of a
monomial ideal $M$. Then we have $\{\rad(M_i) \mid 1 \leq i \leq r\} =
\Ass(S/M)$.  Note that distinct irreducible components may have the same
radical.  Bayer, Peeva and Sturmfels \cite[\S 3]{BPS} give a method for
computing the irreducible decomposition of a generic monomial ideal (in
the old definition).  The generalization of this method by Miller
\cite[Theorem~5.12]{Mil} shows that \cite[Theorem~3.7]{BPS} remains valid
here, as we will show in Theorem 2.2 below.

Recall that $\codim (I) \leq \codim (P) \leq \pdim (S/I) \leq n$ 
for any graded ideal $I \subset S$ and any associated prime $P \in \Ass(S/I)$,
and $\codim (I) = \pdim (S/I)$ if and only if $S/I$ is Cohen-Macaulay. 
There always exists a minimal prime $P \in \Ass(S/I)$ with $\codim (P) = 
\codim (I)$. But in general
there is no $P \in \Ass(S/I)$ with $\codim (P) = \pdim (S/I)$.
For example, if $I = \<x_1,x_2 \> \cap \<x_3, x_4 \>$, then 
$\pdim (S/I) = 3$. 

\begin{thm}\label{gen}
Let $M \subset S$ be a generic monomial ideal.  Then

(a) For each integer $i$ with $\codim (M) < i \leq \pdim (S/M)$, 
there is an {\it embedded} associated prime $P \in
\Ass (S/M)$ with $\codim (P) =i$.

(b) For all $P \in \operatorname{Ass} (S/M)$ there is a chain of
associated primes $P = P_0 \supset P_1 \supset \! \cdots \! \supset P_t$ 
with  $\codim (P_{i}) = \codim (P_{i-1})-1$ for all $i\,$ and
$P_t$ is a minimal prime of $M$.
\end{thm}

\begin{proof}
(a) This was proved by Yanagawa \cite{Y} under the old definition
of genericity.  Using Theorem~\ref{new-def} and \cite[Theorem~5.12]{Mil},
the argument in \cite{Y} also works here.

(b) It suffices to show that for any embedded prime $P$ of $M$ there is
an associated prime $P' \in \Ass(S/M)$ with $\codim (P') = \codim(P) -1$
and $P' \subset P$.  The localization $P_{(P)} $ of $P$ is a maximal
ideal of $S_{(P)}$, and an embedded prime of $M_{(P)}$, so there is a
prime $P'_{(P)} \subset S_{(P)}$ such that $P'_{(P)} \in \Ass
(S/M)_{(P)}$, $\codim (P'_{(P)}) = \codim(P_{(P)}) - 1$ and $P'_{(P)}
\subset P_{(P)}$ by (a) applied to the generic ideal $M_{(P)}$. The
preimage $P' \subset S$ of $P'_{(P)} \subset S_{(P)}$ has the expected
properties.
\end{proof}

\begin{rem} 
Let $M \subset S$ be a generic monomial ideal, and
$P, P' \in \Ass(S/M)$ such that $P \supset P'$ and
$\codim P \geq \codim P' + 2$.  Theorem~\ref{gen} does {\it not} state
that there is an associated prime between $P$ and $P'$. For example, set
$M = \<ac, bd, a^3b^2, a^2b^3\>.$ Then $\<a,b\>, \<a,b,c,d\> \in \Ass
(S/M)$, but there is no associated prime between them.
\end{rem}

Following \cite[\S 3]{BPS}, we next
define the {\it extended Scarf complex} $\Delta_{M^*}$ of $M$. Let
\begin{equation} \label{M*}
  M^* \,\,\, := \,\,\, M + \<x_1^D, \ldots, x_n^D\>
\end{equation}
with $D$ larger than any exponent on any minimal generator of $M$.
We index the new monomials $x_s^D$ just by their variables $x_s$;
so the vertex set of $\Delta_{M^*}$ is a subset of 
$\{1, \ldots, r\} \cup \{x_1,\ldots,x_n \}$. 
This subset is proper if $M$ contains a power of a variable.  Recall
(\cite[Corollary~5.5]{BPS} for the old genericity or
\cite[Proposition~5.16]{Mil} for the new) that $\Delta_{M^*}$ is a
regular triangulation of an $(n-1)$-simplex $\Delta$. The vertex set of
$\Delta$ equals $\{x_1,\ldots,x_n\}$ unless $M$ contains a power of a
variable.  The restriction of $\Delta_{M^*}$ to $\{1, \ldots, r \}$
equals the Scarf complex $\Delta_M$ of $M$.  We next determine the
restriction of $\Delta_{M^*}$ to $\{x_1,\ldots,x_n\}$.

The radical $\rad(M)$ of $M$ is a square-free monomial ideal. 
Let $V(M)$ denote the corresponding 
{\it Stanley-Reisner complex}, which consists of 
all subsets of $\{x_1,\ldots,x_n\}$ which are not 
support sets of  monomials in $M$. Then we have the following:

\begin{lem} \label{rad}
For a generic monomial ideal $M$, the restriction
of the extended Scarf complex $\Delta_{M^*}$  to $\{x_1,\ldots,x_n\}$ 
coincides with the Stanley-Reisner complex $V(M)$. 
\end{lem}

\begin{proof}
Every facet $\sigma$ of $\Delta_{M^*}$ gives an irreducible component of
$M$; see \cite[Theorem~3.7]{BPS} and \cite[Theorem~5.12]{Mil}.  The
radical of that component represents the face $\, \sigma \,\cap\,
\{x_1,\ldots,x_n\}\,$ of $V(M)$. The facets of $V(M)$ arise in this way
from the irreducible components whose associated primes are minimal.
\end{proof}

The following theorem generalizes a result of Yanagawa
\cite[Corollary~2.4]{Y}.  For the definition of shellability, see
\cite[\S III.2]{Sta} or \cite[Lecture 8]{Zie}.

\begin{thm} \label{CMness}
Let $M$ be a generic monomial ideal. If $M$ has no embedded associated
primes, then $M$ is Cohen-Macaulay. In this case, both $\Delta_M$ and
$V(M)$ are shellable. 
\end{thm}

\begin{proof}
The first statement immediately follows from Theorem~\ref{gen}.  For the
second statement we note that all facets $\sigma$ of $\Delta_{M^*}$ have
the following property:
\begin{equation} \label{resttt}
|\sigma \cap \{1, \ldots, r \}| \,=\, \codim M 
\quad \hbox{and} \quad
|\sigma \cap \{x_1, \ldots, x_n \}| \,=\, \dim S/M  .
\end{equation}
In particular, both cardinalities in
(\ref{resttt}) are independent of the facet $\sigma$.
On the other hand,
$\Delta_{M^*}$ is shellable since it is  a regular triangulation of
a simplex. A theorem of Bj\"orner \cite[Theorem~11.13]{Bj}
implies that the restrictions of $\Delta_{M^*}$ to
$\{1,2,\ldots,r\}$ and to $\{x_1,\ldots,x_n\}$ are
both shellable. We are done in view of
Lemma~\ref{rad}.
\end{proof}

\begin{rem} 
\label{justaremark}
(a) The shellability of $\Delta_{M^*}$ also implies the
following result. If $M$ is generic and
$P,P' \in \Ass(S/M)$, then there is a sequence of associated
primes $P = P_0, P_1, \ldots ,$ $ P_t = P'$ with $\codim (P_i +
P_{i-1}) =\min \{ \codim (P_i), \codim (P_{i-1}) \}+1$ for all $1 \leq
i \leq t$. If $M$ is pure dimensional, this simply says that $S/M$ is
connected in codimension $1$.

(b) A shelling of the boundary complex of a polytope can start from a
shelling of the subcomplex consisting of all facets containing a given
face; see \cite[Theorem~8.12]{Zie}.  The complex $V(M)$ of a generic
Cohen-Macaulay monomial ideal $M$ inherits this property, so $V(M)$ has
stronger properties than general shellable complexes.
\end{rem}

Theorem~\ref{CMness} and Remark~\ref{justaremark}
suggest the following combinatorial problems:

\begin{prob}\label{characterize}
(i) Characterize all collections ${\mathcal A}$ of monomial primes for which
there exists a generic monomial ideal $M$ with ${\mathcal A} = \Ass(S/M)$.

(ii) Characterize the Stanley-Reisner complexes $V(M)$
of Cohen-Macaulay generic monomial ideals $M$.
\end{prob}

A necessary condition for (i) is that ${\mathcal A}$
satisfy the connectivity in Remark~\ref{justaremark} (a). 
But this is not sufficient: for instance, take ${\mathcal A}$ to be
the minimal primes of a Stanley-Reisner ring which is
Cohen-Macaulay but whose simplicial complex not shellable.

For the problem (ii), the Cohen-Macaulayness assumption is essential.
Since for all simplicial complex $\Sigma \subset 2^n$, there is a (not
necessarily Cohen-Macaulay) generic monomial ideal $M$ such that $V(M) =
\Sigma$.  By Theorem~\ref{CMness}, shellability is a necessary condition
for the problem (ii), but it is not sufficient as
Remark~\ref{justaremark} (b) shows.

If we put further restrictions on the generators of a generic monomial
ideal $M$, then, since the extended Scarf complex $\Delta_{M^*}$ is a
triangulation of a simplex, we can apply Stanley's theory of local
$h$-vectors \cite{Sta}.  The next two results will
be reinterpreted in Section~\ref{cogeneric} in terms of cogeneric ideals
using Alexander duality \cite{Mil}.

Again let $M^*$ be as in \eqref{M*}, and define the {\it excess} of a
face $\sigma \in \Delta_{M^*}$ to be $e(\sigma) := \#\supp(m_\sigma) -
\#\sigma$.  This agrees, in our situation, with the definition of excess
in \cite{Sta}.

\begin{thm}\label{size=c}
If $M$ is generic and all $r$ generators $m_1,\ldots,m_r$ have support of
size $c$, i.e.\ $\# \supp(m_i) = c$ for all $i$, then $M$ has at least
$(c-1) \cdot r + 1$ irreducible components.
\end{thm}

\begin{exmp} \label{false}
This is false without the assumption that $M$ is generic.  For instance,
the non-generic monomial ideal $M = \<x_1,y_1\> \cap \ldots \cap
\<x_n,y_n\>$ has $r=2^n$ generators, and each generator has support of
size $c = n$, but $M$ has only $n$ irreducible components.
\end{exmp}

\begin{proof} 
If $c =1$, there is nothing to prove, so we may assume that $c \geq 2$. 
Set $\Gamma = \Delta_{M^*}$. The 
hypothesis on the generators of $M$ means that $\Gamma$ has $n$ vertices of 
excess $0$ and $r$ vertices of excess $c-1$.  To prove the assertion, we 
use the decomposition 
\begin{equation} \label{h=l'} 
  h(\Gamma,x) \ = \ \sum_{W \in \Delta} \ell_W(\Gamma_W,x)
\end{equation}   
of the $h$-polynomial of $\Gamma$ into local $h$-polynomials
\cite[eqn.~(3)]{Sta}.
Here $\Delta$ denotes the simplex on $\{x_1,\ldots,x_n\}$
and $\Gamma_W$ the restriction of $\Gamma$ to a face $W$ of $\Delta$. We have
\begin{equation} \label{l=1'}
  \ell_W(\Gamma_W,x) =  1  \qquad \text{if $W = \emptyset$.}
\end{equation}
Next, we consider the case $\# W = c$.  In the $\Gamma_W$, the vertices
corresponding to generators of $M$ have excess $c-1$, and all other faces
have excess less than $c-1$.  So we have
\begin{equation} \label{lexpand'}
  \ell_W(\Gamma_W,x) =  \ell_1(\Gamma_W)x + \ell_2(\Gamma_W) x^2 + \cdots + 
  \ell_{c-1}(\Gamma_W) x^{c-1} \qquad \text{if $\#W = c$,}
\end{equation} 
where $\ell_1(\Gamma_W)$ is the number of generators of $M$ whose support 
corresponds to the face $W$ of $\Delta$ by \cite[Example~2.3(f)]{Sta}.  
Moreover $\ell_i(\Gamma_W) \geq \ell_1 (\Gamma_W)$ for all $1 \leq i \leq c-1$ 
by \cite[Theorem~5.2 and Theorem~3.3]{Sta}.

The coefficients of $\ell_W(\Gamma_W,x)$ are non-negative for all $W \in
\Delta$ by \cite[Corollary~4.7]{Sta}. We now substitute the expressions in
\eqref{l=1'} and \eqref{lexpand'} into the sum on the right hand side of
\eqref{h=l'}, and then we evaluate at $x=1$.  The number of irreducible
components of $M$ equals the number $\, f_{n-1}(\Gamma) \, = \, h(\Gamma,
1) \,$ of facets of $\Gamma$ by \cite[Theorem~5.12]{Mil}, hence
$$
  h(\Gamma,1)
  \,\,\, \geq \,\,\,
	  1 + \sum_{\#W = c}(\sum_{i=1}^{c-1}\ell_i(\Gamma_W))
  \,\,\, \geq \,\,\,
	1 + \sum_{\#W = c} (c-1) \cdot \ell_1 (\Gamma_W) 
  \,\,\,  =  \,\,\, 
	(c-1) \cdot r + 1.
$$
This yields the desired inequality.
\end{proof}

The inequality in Theorem~\ref{size=c} is sharp for all $c$ 
and $r$; see Example~\ref{optimal} below.

\begin{prop} \label{size=2}
Let $M$ be a generic monomial ideal with $r$ generators each of which is 
a bivariate monomial.  Then $M$ has exactly $r + 1$ irreducible components 
if and only if $\#\supp(m_\sigma) \leq 3$ for all edges $\sigma \in
\Delta_M$.
\end{prop}

\begin{proof} 
By the assumption, $\Delta_{M^*}$ has $n$ vertices of excess 0 and $r$ 
vertices of excess 1.  Adding a vertex to any face of $\Delta_{M^*}$ increases 
the excess by at most $1$, so we conclude that the equality 
$\{\sigma \in \Delta_{M^*} \mid \#\sigma = e(\sigma)\} = \{\,\emptyset, 
\{1\}, \{2\}, \cdots, \{r\} \, \}$ holds if and only if each edge of 
$\Delta_M$ has excess at most 1, equivalently, support of size at most $3$. 
The result is now an immediate consequence of \cite[Proposition~3.4]{Sta}.
\end{proof}

\section{Initial Ideals of Lattice Ideals}

One motivation for our new definition of genericity for monomial ideals
is consistency with the notion of genericity for lattice ideals
introduced in \cite{PS2}.  It is the purpose of this section to establish
this connection.  We fix a sublattice ${\LL}$ of ${\mathbb Z}^n$ which
contains no nonnegative vectors.  The {\it lattice ideal} $I_{\LL}$
associated to ${\LL}$ is defined by
$$  I_{\LL} \,\,\, := \,\,\, \< \,{\bf x}^{{\bf a}}
 - {\bf x}^{{\bf b}} \,| \ {\bf
  a},{\bf b}\in {\mathbb N}^n\ \ \mbox{and}\ \ {\bf a} - {\bf b}\in {\LL}
  \, \> \quad  \subset  \,\,\,\, S, $$
where ${\bf x}^{\bf a} = x_1^{a_1} \cdots x_n^{a_n}$ for ${\bf a} =
(a_1,\ldots,a_n) \in {\mathbb N}^n$. The ideal $I_{\LL}$ is homogeneous
with respect to some grading where $\deg(x_s)$ is a positive integer for
each $s$.  We have $\codim(I_{\LL}) = \rank (\LL)$.  Moreover, the
ring $S/I_{\LL}$ has a fine grading by ${\mathbb Z}^n/{\LL}\,$
(cf.~\cite{PS}).

The following three conditions are equivalent: (a) The abelian group
 ${\mathbb Z}^n/{\LL}$ is torsion free, (b) $I_{\LL}$ is a prime
ideal, and (c) $I_{\LL}$ is a toric ideal (i.e., $S/I_{\LL}$ is an affine
semigroup ring).  Even if $I_{\LL}$ is not prime, all monomials
are non-zero divisors of $S/I_{\LL}$, and all associated primes of
$I_{\LL}$ have the same codimension.  If $I_A$ is the toric ideal
of an integer matrix $A$, as defined in \cite{Stu},
then $I_A $ coincides with the lattice ideal $I_{\LL}$ 
where $\LL  \subset {\mathbb Z}^n$
is the  kernel of $A$.

Following Peeva and Sturmfels \cite{PS2}, we call a
lattice ideal $I_{\LL}$ 
{\it generic} if it is generated by binomials with full support, i.e.,
$$ \,  I_{\LL}  \,\,\,=  \,\,\,  \< \, 
	{\bf x}^{{\bf a}_1} -  {\bf x}^{{\bf b}_1}, \,
	{\bf x}^{{\bf a}_2} -  {\bf x}^{{\bf b}_2}, \, \ldots \,, \,
	{\bf x}^{{\bf a}_r} -  {\bf x}^{{\bf b}_r} \, \> \,$$
where none of the $r$ vectors ${\bf a}_i - {\bf b}_i \in {\mathbb Z}^n$
has a zero coordinate.

\begin{thm} \label{ingen}
Let $I_{\LL}$ be a generic lattice ideal, and  $M$
the initial ideal of $I_{\LL}$ with respect to a 
reverse lexicographic term order.  Then $M$ is a generic monomial ideal.  
\end{thm}

\begin{proof} 
Set $M = \init_{revlex}(I_{\LL}) =  \<m_1, \ldots, m_r\>$.
Gasharov, Peeva and Welker \cite{GPW} proved that
the algebraic Scarf complex $F_{\Delta_M}$ is a minimal
free resolution of $S/M$. Using Theorem~\ref{new-def},
 it suffices to prove that no variable $x_s$
appears with the same non-zero exponent in $m_i$ and $m_j$ for any $i \ne
j$ with $\{ i,j \} \in \Delta_M$.  Assume the
contrary, that is, $\deg_{x_s} m_i = \deg_{x_s} m_j > 0$ for some
$\{i,j\} \in \Delta_M$. By \cite[Theorem 5.2]{PS2}, 
there are three monomials $m_i',
m_j', m_l' \in S$ satisfying the following conditions.

(a) $\{ m_i', m_j', m_l' \}$ is a  {\it basic fiber} (see \cite[\S 2]{PS2}), 
in particular, $\gcd (m_i', m_j', m_l') = 1$. 

(b) $m_i = \frac{m_i'}{\gcd (m_i', m_l')}\,\,$ and $\,\,m_j = \frac{m_j'}{\gcd
(m_j', m_l')}$.

By (b), we have $\deg_{x_s} (m_i') \geq \deg_{x_s} (m_i) > 0$ and
$\deg_{x_s} (m_j') \geq \deg_{x_s} (m_j) > 0$.  Since $\gcd (m_i', m_j',
m_l') = 1$, we have $\deg_{x_s} m_l' = 0$.  So $\deg_{x_s} m_i' =
\deg_{x_s} m_i = \deg_{x_s} m_j = \deg_{x_s} m_j'$.  
Combining property (a) with \cite[Theorem 3.2]{PS2},
we see that the binomial
$$
  \frac{m_i'}{\gcd (m_i', m_j')} - \frac{m_j'} {\gcd (m_i', m_j')}
$$
is a minimal generator of $I_{\LL}$.  Since $\deg_{x_s} m_i' = \deg_{x_s}
m_j'$, the variable $x_s$ does not appear in the above binomial. This
contradicts the genericity of $I_{\LL}$.
\end{proof}

\begin{exmp}\label{gen-lat-ex}
Theorem~\ref{ingen} is false for the
old definition of ``generic monomial ideal'' given in  
\cite{BPS}. For example, consider the following generic lattice ideal
in $k[a,b,c, d]$:
$$  I_\LL \,\,= \,\, \<a^4 - b c d, 
  a^3 c^2 - b^2 d^2,  a^2 b^3 - c^2 d^2, a b^2 c - d^3, 
  b^4 - a^2cd,\, b^3 c^2 - a^3 d^2, 
  c^3 - a b d \> $$ 
This ideal was featured in  \cite[Example 4.5]{PS2};
it defines the toric curve
$ (t^{20},t^{24},t^{25},t^{31})$.  Consider a 
reverse lexicographic term order with $a > b > c > d$.  Then $\,M =
\<a^4, a^3 c^2, a^2 b^3, a b^2 c, b^4, b^3 c^2, c^3 \>$.  Since 
$a^3c^2$ and $b^3c^2$ are minimal
generators of $M$, it is not generic in the sense of \cite{BPS}.  
But $M$ satisfies Definition \ref{gendef} since $a b^2 c \in M$.  \qed
\end{exmp}

\vskip .1cm

An important problem in combinatorial commutative algebra is to
characterize those monomial ideals which are initial ideals of lattice
ideals.  The recent ``Chain Theorem'' of Ho\c{s}ten and Thomas \cite{HT}
provides a remarkable necessary condition.

\begin{thm}[Ho\c sten--Thomas \cite{HT}] \label{chain}
Let $M$ be the initial ideal of a lattice ideal $I_{\LL}$ with respect to
any term order.  For each $P \in \operatorname{Ass} (S/M)$, there is a
chain of associated primes $P = P_0 \supset P_1 \supset \cdots \supset
P_t$ of $M$ such that $P_t$ is a minimal prime and $\codim (P_{i}) =
\codim (P_{i-1})-1$ for all $i$.
\end{thm}

In other words, initial ideals of lattice ideals satisfy conclusion (b)
of Theorem~\ref{gen}, even if they are not generic. We do not know
whether part (a) holds as well.

\begin{conj}\label{prob}
Let $M$ be the initial ideal of $I_{\LL}$ with respect to some term
order. Then there is an associated prime $P \in \Ass (S/M)$ with $\codim
(P) = \pdim (S/M)$.
\end{conj}

\begin{cor} Conjecture~\ref{prob} holds for the reverse lexicographic 
term order if the lattice ideal $I_{\LL}$ is generic.
\end{cor}

\begin{proof}
Immediate from Theorem~\ref{gen} and Theorem~\ref{ingen}.
\end{proof}

The following result appears implicitly in
the work of Ho\c{s}ten-Thomas \cite{HT} and Peeva-Sturmfels \cite{PS}.

\begin{lem}\label{pdim}
Let $M$ be the initial ideal of a lattice ideal $I_{\LL}$
with respect to any term order.  Then we have 
$\pdim(S/M) \leq 2^c-1$ where $c := \codim I_{\LL} = \codim M.$
\end{lem}

\begin{proof}
Following \cite[Algorithm 8.2]{PS}, we construct a lattice ideal $I_{\LL'} $ 
in $\, S[t] = k[x_1,\ldots,x_n,t]$ whose images under the
substitutions $t=1$ and $t=0$ are $I_{\LL}$ and $M$ respectively. 
Moreover $t$ is a non-zero divisor of $S[t]/I_{\LL'}$, and  
the codimension of $I_{\LL'}$ in $S[t]$ is equal to $\codim (I_{\LL})$.  
Since $S/M = S[t]/(I_{\LL'} + \<t\>)$, we have 
$\pdim(S/M) = \operatorname{proj-dim}_{S[t]} (S[t]/I_{\LL'}) \leq 
2^c - 1$. The last inequality follows from \cite[Theorem~2.3]{PS}. 
\end{proof}

We note that Conjecture~\ref{prob} is also true in codimension 2:  

\begin{prop} Conjecture \ref{prob} holds for any term order if 
$\codim(I_{\LL})=2$. 
\end{prop}

\begin{proof} 
By Lemma~\ref{pdim}, $\pdim (S/M) \leq 3$. 
We may assume  $\pdim (S/M) = 3$, because otherwise $M$ is Cohen-Macaulay 
and there is nothing to prove. Then there exists a syzygy quadrangle as 
in \cite[\S 3]{PS} for the planar configuration of $n+1$ vectors representing 
the ideal $I_{\LL'}$ from Lemma~\ref{pdim}. 
This quadrangle defines a lattice point free polytope as in 
\cite[\S 2]{HT}, and from the explicit primary decomposition given 
by Ho\c{s}ten  and Thomas \cite[Theorem 4.2]{HT} we see that $M$ has 
an associated  prime of codimension $3$. 
\end{proof}

For an ideal $I \subset S$, it is well-known that 
$\pdim (S/I) \leq \pdim(S/\init(I))$. 
This inequality can be strict even in the codimension 2 toric ideal 
case. Set $I_{\LL} := \< ac-b^2, ad-bc, bd-c^2 \> \subset S = k[a,b,c,d]$ 
be the  defining ideal of the twisted cubic curve in ${\mathbb P}^3$. 
$S/I_{\LL}$ is normal and Cohen-Macaulay. 
The ideal $I_{\LL}$ has eight 
distinct initial ideals, when we consider all possible term orders 
(see \cite[\S 4]{Stu0}). Four of them are not Cohen-Macaulay and have 
embedded associated primes of codimension 3.
\begin{rem}

Let $M \subset S$ be a {\it Borel fixed} monomial ideal 
(cf. \cite[\S 15.9]{E}).  In general, Borel fixed ideals are far from generic. 
But it is easy to see that there is an associated prime $P \in \Ass(S/M)$ 
with $\codim (P) = \pdim (S/M)$. Hence a Borel fixed ideal $M$ satisfies 
the conclusion of Conjecture~\ref{prob}. Therefore the {\it generic initial 
ideal} (cf. \cite[\S 15]{E}) of a homogeneous ideal $I \subset S$ satisfies 
the conclusion of the conjecture, when $\operatorname{char} k = 0$. 
But Borel fixed ideals  may fail to satisfy the conclusion of 
Theorem~\ref{chain}. For instance, take $M = \<x^2,xy,xz\> = \<x\> \cap 
\<x^2,y,z\>$.
\end{rem}

\section{A Study of Cogeneric Monomial Ideals} \label{cogeneric}

Cogeneric monomial ideals were introduced in Definition~\ref{cogendef}.
As  with genericity, our definition of 
cogenericity is slightly different from the original one of \cite{Stu2}. 
In Theorem~\ref{coscarf} we shall see that 
the result of \cite{Stu2}, an explicit description of the
minimal free resolution of a cogeneric monomial ideal, is still true here.
In fact, Alexander duality for arbitrary monomial ideals \cite{Mil} allows 
us to shorten the construction of this resolution and clarify its 
relation to Theorem  \ref{new-def}. 
For the reader's convenience, we briefly recall the 
definitions pertaining to Alexander duality.  For details see \cite{Mil}.

The maximal $\NN^n$-graded ideal $\<x_1, \ldots, x_n\> \subset S$ will be
denoted by $\mm$.  Monomials and irreducible monomial ideals may each be
specified by a single vector $\bb = (b_1,\ldots,b_n) \in \NN^n$, so we
will write ${\bf x}^\bb = x_1^{b_1} \cdots x_n^{b_n}$ and $\mm^\bb =
\<x_s^{b_s} \mid b_s \geq 1\>$.  Given a vector $\aa = (a_1,\ldots,a_n)$
such that $b_s \leq a_s$ for all $s$, we define the Alexander dual vector
$\bb^\aa$ with respect to $\aa$ by setting its $s^{\,\rm th}$ coordinate
to be
$$
  (\bb^\aa)_s = \left\{
	\begin{array}{ll}
	a_s + 1 - b_s	& {\rm if}\ b_s \geq 1 \cr
	0		& {\rm if}\ b_s = 0.
	\end{array}
  \right.
$$
Whenever we deal with Alexander duality, we assume that we are given a
vector $\aa$ such that for each $s$, the integer $a_s$ is 
larger than or equal to 
the $s^{\,\rm th}$ coordinate of any minimal monomial generator of $M$.
This implies that $a_s$ is also larger than or equal to the $s^{\,\rm th}$ 
coordinate of any irreducible component of $M$, and vice versa.
The {\it Alexander dual ideal}
$M^\aa$ of $M$ with respect to $\aa$ is defined by
\begin{eqnarray*}
  M^\aa & = &
	\<\xx^{\bb^{\scriptstyle \aa}} \, | \, \text{$\mm^\bb$ is an
	irreducible component of $M$} \>\\
  & = &
	\bigcap \{\mm^{\cc^{\scriptstyle \aa}} \, | \, \text{$\xx^\cc$ is
	a minimal generator of $M$} \}. 
\end{eqnarray*}
That these two formulas give the same ideal is not obvious;
it is equivalent to  $(M^\aa)^\aa = M$.  
It follows from these statements that $M$ is generic 
if and only if $M^\aa$ is cogeneric.

\begin{exmp}
The following monomial ideal in $S = k[x,y,z]$ is cogeneric:
$$
  M \,\,= \,\, \<yz^2, xz^2, y^2z, xy^2, x^2\> \,\,= \,\, \<x,y\> \cap
  \<x^2,y^2,z^2\> \cap \<x,z\>.
$$
Its Alexander dual with respect to $\aa = (2,2,2)$ is generic:
$$
  M^\aa \,\,= \,\,\<x^2y^2, xyz, x^2z^2\>\,\, = \,\, \<y^2,z\> \cap
  \<x^2,z\> \cap \<y,z^2\> \cap \<x^2,y\> \cap \<x\>.
$$

\end{exmp}

\begin{exmp}[{\cite[Examples~1.9, 5.22]{Mil}}] \label{permuu}
If $M$ is the tree ideal of Example~\ref{tree} and $\aa =
(n,\ldots,n)$, then its Alexander dual $M^\aa$ is the {\it permutahedron
ideal}:
$$
  M^\aa \quad = \quad
  \< \,
   x_1^{\pi(1)}
   x_2^{\pi(2)} \cdots
   x_n^{\pi(n)} \,: \,
  \pi \,\, \text{is a permutation of} \,\,
  \{1,2,\ldots,n\} \, \>.
$$
Thus the permutahedron ideal is cogeneric.  Its minimal free resolution
is the {\it hull resolution}, which is cellular and supported on a
permutahedron \cite[Example~1.9]{BS}.  The following discussion
reinterprets this resolution as a co-Scarf complex. \qed
\end{exmp}

\begin{dfn}\label{co-Scarf}
Let $M=\bigcap_{i=1}^r M_i$ be a cogeneric monomial ideal. 
Set $\aa = (D - 1, \ldots, D - 1)$ with $D$ larger than any exponent on any 
minimal generator of $M$. The Alexander dual ideal $M^\aa$ is minimally 
generated by monomials $m_1, \ldots, m_r$, where $m_i = \xx^{{\bb_i}^
{\scriptstyle \aa}}$ for $M_i = \mm^{\bb_i}$. We define the 
{\it co-Scarf complex} 
$\Delta^\aa_M$ to be the extended Scarf complex of $M^\aa$. 
More precisely, we set $(M^\aa)^* := M^\aa + \<x_1^D, \ldots, x_n^D\>$ and  
$\Delta^\aa_M$ the Scarf complex of $(M^\aa)^*$. 
Since we index a new monomial $x_s^D$ just by $x_s$, we see that
$\Delta^\aa_M$ is a simplicial complex 
on (a subset of) $\,\{1, \ldots, r, x_1,\ldots,x_n \}$.
\end{dfn}

\begin{rem}
(a) There is nothing special about our choice of $\aa$, except that it
makes for convenient notation.  Everything we do with $\Delta_M^\aa$ is
independent of which sufficiently large $\aa$ is chosen.  In particular,
the regular triangulation of the $(n-1)$-simplex is independent of $\aa$,
as is the algebraic co-Scarf complex (Definition~\ref{alg}) it
determines.  We therefore set $\aa = (D-1,\ldots,D-1)$ for the remainder
of this section.

(b) For $\sigma \subseteq \{1,\ldots,r\}$, let $M_{\sigma}$ be the
irreducible monomial ideal $\sum_{i \in \sigma} M_i$.  
Then $m_\sigma = \xx^{\bb^{\scriptstyle \aa}}$ if $M_\sigma = \mm^\bb$,
and $\Delta_M^\aa \cap \{1, \ldots, r\}= \{\sigma \subset \{1, \ldots, r\}
\mid M_{\tau} \ne M_{\sigma} \ \text{for all} \ \tau \ne \sigma \}$ is
just the Scarf complex of $M^\aa$.
\end{rem}

A face $\sigma$ of the co-Scarf complex $\Delta_M^\aa$ fails to be in the 
(topological) boundary $\partial\Delta_M^\aa$ of $\Delta_M^\aa$ if and only if 
the monomial $m_\sigma$ has full support, where 
$m_{\sigma}$ is $\lcm ( m_i \mid i \in \sigma )$ under the notation of 
Definition~\ref{co-Scarf}. 
Such a face will be called an {\it interior face} of $\Delta_M^\aa$.  
The set ${\rm int}(\Delta_M^\aa)$
of interior faces is closed under taking supersets; that is, ${\rm
int}(\Delta_M^\aa)$ is a {\it simplicial cocomplex}.  Just as the
algebraic Scarf complex is constructed from $\Delta_M$ for generic $M$,
we construct an algebraic free complex from ${\rm int}(\Delta_M^\aa)$,
but this time we use the coboundary map instead of the boundary map.  The
following is a special kind of {\it relative cocellular resolution} (in
fact a cohull resolution) \cite[\S5]{Mil}.

\begin{dfn} \label{alg}
Let $\DD = (D,\ldots,D) \in \NN^n$ and $S(\aa_\sigma-\DD)$ be the free
$S$-module with one generator $e^*_\sigma$ in multidegree $\DD -
\aa_\sigma$.  The {\it algebraic co-Scarf complex} $F^{\Delta_M^\aa}$ of
$M$ is the free $S$-module
$$
  \bigoplus_{\sigma \in {\rm int}(\Delta_M^\aa) } \!\! S(\aa_\sigma - \DD)
  \quad \text{with differential} \quad 
  d^*(e^*_\sigma)\, = \!\!\!\!\!\! \sum_{\substack{i \not\in \sigma \\ \sigma
  \cup \{i\} \in {\rm int}(\Delta_M^\aa)}} \!\!\!\!\!\!
  \operatorname{sign}(i,\sigma \cup \{i\}) \cdot \frac{m_{\sigma \cup
  \{i\}}}{m_\sigma} \cdot e^*_{\sigma \cup \{i\}}
$$
where $\operatorname{sign}(i,\sigma \cup \{i\})$ is $(-1)^{j+1}$ if $i$
is the $j$-th element in the ordering of $\sigma \cup \{i\}$.  Put the
summand $S(\aa_\sigma - \DD)$ in homological degree $n - \#\sigma = n -
\dim(\sigma) - 1$.
\end{dfn}

\begin{thm} \label{coscarf}
If $M$ is a cogeneric monomial ideal, then the algebraic co-Scarf
 complex $F^{\Delta_M^\aa}$ equals the minimal free resolution of $M$
 over $S$.  In particular, $M$ is minimally generated by the set of
 monomials $\,\{\xx^{\DD - \aa_{\scriptstyle \sigma}} \mid \sigma$ is
 a facet of $\Delta_M^\aa\}$.
\end{thm}

\begin{proof}
This  follows from Proposition~\ref{hull} and
\cite[Theorem~5.8]{Mil}.
\end{proof}

\noindent {\bf Example 4.1} {\sl (continued) }
For the cogeneric ideal $M = \<x,y\> \cap \<x^2,y^2,z^2\> \cap
\<x,z\>$, the interior faces of $\Delta_M^\aa$ are $\{2\}$, $\{1,2\}$,
$\{2,3\}$, $\{2,x\}$, $\{2,y\}$, $\{2,z\}$, $\{1,2,x\}$, $\{1,2,y\}$,
$\{2,3,x\}$, $\{2,3,z\}$ and $\{2,y,z\}$.  The co-Scarf resolution 
is $0 \to S \to S^5 \to S^5 \to M \to 0$.  The generators of $M$
have exponent vectors $\,\DD -
\aa_{\{1,2,x\}} = (0,1,2)$, $\DD - \aa_{\{1,2,y\}} = (1,0,2)$, $\DD -
\aa_{\{2,3,x\}} = (0,2,1)$, $\DD - \aa_{\{2,3,z\}} = (1,2,0)$ and $\DD -
\aa_{\{2,y,z\}} = (2,0,0)$. 

\vskip .3cm

We saw in Theorem~\ref{CMness} that for generic monomial ideals, the 
Cohen-Macaulay condition is equivalent to the much weaker condition of 
purity (all associated primes have the same dimension).  For cogeneric 
monomial ideals, on the other hand, purity is obviously too easy to 
attain.  Nonetheless, a cogeneric ideal is forced to be Cohen-Macaulay by 
{\it a priori} much weaker conditions.  Before stating these in 
Theorem~\ref{CMness2}, we characterize depth for cogeneric ideals using a 
polyhedral criterion. 

\begin{lem} \label{depth} 
Let $M$ be a cogeneric monomial ideal.  Then $\depth(S/M) \leq d$ if and 
only if the co-Scarf complex $\Delta_M^\aa$ has an interior face of 
dimension $d$. 
\end{lem}

\begin{proof}
By Theorem~\ref{coscarf}, the shifted augmentation $F^{\Delta_M^\aa} \to S$ 
(obtained by including $\operatorname{coker}(F^{\Delta_M^\aa}) = M$ into 
$S$ and shifting homological degrees up one) is a minimal free resolution 
of $S/M$.  The co-Scarf complex $\Delta_M^\aa$ has an interior face of 
dimension $d$ if and only if this shifted augmented complex is nonzero in 
homological dimension $n-d$.  The lemma now follows from the 
Auslander-Buchsbaum formula. 
\end{proof}

Recall that a module $N$ satisfies Serre's condition $(S_k)$ if for every
prime $P \subset S$, $\depth(N_P) < k$ $\Rightarrow$ $\depth(N_P) =
\dim(N_P)$.  Using \cite[Chapter~2.1]{BH} and homogeneous localization,
it follows that if $S/M$ satisfies $(S_k)$ then
\begin{equation} \label{Sk}
  \depth((S/M)_{(P)}) < k \quad \Longrightarrow \quad
  \dim((S/M)_{(P)}) = \depth((S/M)_{(P)}).
\end{equation}
Observe that $M_{(P)}$ is cogeneric if $M$ is, in analogy to Remark 2.1.
For condition (d) below,
recall the definition of {\it excess} from before Theorem~\ref{size=c}.

\begin{thm} \label{CMness2}
Let $M \subset S$ be a cogeneric monomial ideal of codimension $c$ with the
irreducible decomposition $M = \bigcap_{i=1}^{r} M_i$.  Then the
following conditions are equivalent.

(a) $S/M$ is Cohen-Macaulay.  

(b) $S/M$ satisfies Serre's condition $(S_2)$. 

(c) $\codim M_i = c$ for all $i$, and $\codim (M_i + M_j) \leq c + 1$
for all edges $\{ i, j \} \in \Delta_{M}^\aa$.

(d) Every face of $\Delta_M^\aa$ has excess $< c$.

(e) $\Delta_M^\aa$ has no interior faces of dimension $< n-c$.
\end{thm}

\begin{proof} (a) $\Rightarrow$ (b) : Cohen-Macaulay $\Leftrightarrow (S_k)$
for all $k$.

(b) $\Rightarrow$ (c) : The initial equality follows from
\cite[Remark~2.4.1]{H}, so it suffices to prove the inequality.  Suppose
$i \ne j$ with $\{i,j\} \in \Delta_{M}^\aa$.  Let $P = \rad (M_i + M_j)$,
and denote by $F$ the face of 
$\,\Delta = 2^{ \{x_1,\ldots,x_n \} } \,$ whose vertices are 
the variables in $P$. By \cite[Proposition~4.6]{Mil}, 
the co-Scarf complex of $M_{(P)}$ is, as a triangulation
of the simplex $2^F$, the restriction $(\Delta_M^\aa)_F$ of the
triangulation $\Delta_M^\aa$ to $2^F$.  By our choice of $F$, $\{i,j\}$ is
an interior edge of $(\Delta_M^\aa)_F$, so Lemma~\ref{depth} implies that
$\depth((S/M)_{(P)}) \leq 1$, whence \eqref{Sk} implies that $\dim
((S/M)_{(P)}) \leq 1$.  Equivalently, $\codim(M_i + M_j) \leq c + 1$.

(c) $\Rightarrow$ (d) : The purity of the irreducible components means
that all vertices have excess $c-1$ or 0, while the condition on the
edges implies that the excess of a nonempty face can only decrease or
remain the same upon the addition of a vertex.

(d) $\Rightarrow$ (e) : In particular, the interior faces have excess less
than $c$.

(e) $\Rightarrow$ (a) : Lemma~\ref{depth}.
\end{proof}

\begin{rem} \label{hart}
(a) Hartshorne~\cite{H} proved that a catenary local ring satisfying
Serre's condition $(S_2)$ is pure and connected in codimension $1$.  The
converse is not true even for cogeneric monomial ideals. If we take $M =
\<x,y^2\> \cap \<y,z\> \cap \<z^2,w\>$ then 
  $S/M$ is pure and connected in codimension
$1$, but does not satisfy the condition $(S_2)$; in fact, $\depth (S/M) =
1$. On the other hand, $M' = \<x,y\> \cap \<y^2,z^2\> \cap \<z,w\>$ is
Cohen-Macaulay, although $\Ass(M) = \Ass(M')$.

(b) Let $I$ be a squarefree monomial ideal and $I^\vee =
I^{(1,\ldots,1)}$ its Alexander dual.
Eagon and Reiner~\cite{ER} proved that $S/I$ is Cohen-Macaulay if and only if 
$S/I^{\vee}$ has linear free resolution.  In \cite{Y2}, it is  proved that 
$S/I$ satisfies the $(S_2)$ condition if and only if all minimal generators 
of $I^{\vee}$ have the same degree and all minimal {\it first} syzygies are 
linear. So the equivalence between 
(b) and (c) of Theorem \ref{CMness2} is quite natural, since an edge 
$\{i,j\} \in \Delta_M^\aa$ corresponds to a first syzygy of $M^\aa$. 
But the $(S_2)$ condition is much weaker than 
Cohen-Macaulayness for squarefree monomial ideals.  
\end{rem}

The above theorem and remark leads to a natural question. 

\begin{prob}
Which Cohen-Macaulay simplicial complexes have Stanley-Reisner ideal
$\,\rad(M)\,$ for some Cohen-Macaulay cogeneric
monomial ideal $M$?
\end{prob}

Recall that the {\it type} of a Cohen-Macaulay quotient $S/M$ is the nonzero 
total Betti number of highest homological degree; if $M$ is cogeneric then 
this Betti number equals the number of interior faces of minimal dimension in
$\Delta_M^\aa$ by Theorem~\ref{coscarf}.

\begin{thm}\label{CMtype}
Let $M$ be a Cohen-Macaulay cogeneric monomial ideal of codimension 
$\geq2$.  The type of $S/M$ is at least the number of irreducible
components of $M$.
\end{thm}

Recall that $S/M$ is {\it Gorenstein} if its Cohen-Macaulay type equals
$1$.  This implies:

\begin{cor}\label{Gor}
Let $M$ be a cogeneric monomial ideal.  Then $S/M$ is Gorenstein if and
only if $M$ is either a principal ideal or an irreducible ideal.
\end{cor}

\begin{rem}
In the generic monomial ideal case, we have the opposite inequality
to the one in Theorem~\ref{CMtype}.
More precisely, if $M$ is Cohen-Macaulay and generic then
\begin{align*}
\text{Cohen-Macaulay type of $S/M$} 
&\quad =\quad \# \{ \text{facets of the Scarf complex $\Delta_M$} \} \\
\leq \quad  \# \{ \text{facets of  $\Delta_{M^*}$} \} 
&\quad = \quad \# \{\text{irreducible components of $M$} \}, 
\end{align*}
because the map $\,\Delta_{M^*} \rightarrow \Delta_M,\, 
\sigma  \mapsto \sigma \cap \{1,\ldots,r\} \,$ is surjective on facets.
Also here, $S/M$ is Gorenstein if and only if it is complete 
intersection \cite[Corollary~2.11]{Y}.
\end{rem}

We present two proofs of Theorem~\ref{CMtype}.  The first is algebraic
and uses Alexander duality, in particular the following result.  For
notation, define $\bb \cdot F \in \NN^n$, for $F \subseteq
\{1,\ldots,n\}$ and $\bb \in \NN^n$, to have $s^{\,\rm th}$ coordinate
$b_s$ if $s \in F$ and 0 otherwise.  Also, set $\beta_{i,\bb}(M) =
\dim_k(\Tor_i^S(M,k))_\bb$, the $i^{\,\rm th}$ Betti number of $M$ in
$\ZZ^n$-degree $\bb$.

\begin{thm}[E. Miller {\cite[Theorem~4.13]{Mil}}] \label{thm:inequality}
Let $M \subset S$ be any monomial ideal and let $F \subseteq \{1, \ldots, n
\}$.  If $\supp(\bb) = F$ and $b_s \leq a_s$ for all $s$, then
$$
  \beta_{i \, , \, \bb^{\aa} \,}(M^\aa) \  \leq \ 
  \sum_{\substack{\cc \in \NN^n \\ \cc \cdot F = \bb}}
  \beta_{\#F - i - 1 \, , \, \cc \,}(M).
$$
\end{thm}

\vskip .5mm
\noindent
{\it Proof of Theorem~\ref{CMtype}.}
Let $\irr(S/M)$ denote the set of vectors $\bb \in \NN^n$
for which $\mm^{\bb} $ is an irreducible component of
$M$.  For any $\cc \in \NN^n$, we define
$$  \gamma_\cc \ :=\ 
  \#\{F \subseteq \{1, \ldots, n\} \, | \, \cc \cdot F \in \irr (S/M) \}.
$$
Set $d = \codim(M)$. The first aim is to show that
\begin{equation} \label{eqn:1}
  \#\irr(S/M) \ \leq \ 
  \sum_{\cc \in \NN^n} \gamma_\cc \cdot \beta_{d-1 \, , \, \cc \, }(M).
\end{equation}
In fact, this inequality holds even if $M$ is not cogeneric: by the
construction of $M^\aa$,
$$
  \#\irr(S/M)
  \ = \ 
	\sum_{\bb \in \irr (S/M)} \beta_{0 \, , \, \bb^{\aa}\,}(M^\aa)
  \ = \ 
	\sum_{\bb \in \NN^n} \beta_{0 \, , \, \bb^{\aa} \,}(M^\aa).
$$ 
Since $S/M$ is Cohen-Macaulay of codimension $d$, each $\bb \in \irr
(S/M)$ has precisely $d$ non-zero coordinates, and
$\beta_{i\,,\,\cc\,}(M) = 0$ for $i \geq d$.  Thus
Theorem~\ref{thm:inequality} specializes to
$$
	\beta_{0\, , \, \bb^{\aa}\,}(M^\aa) \ \leq \
	\sum_{\cc \cdot F \,=\, \bb} \beta_{d - 1\,,\,\cc\,}(M)
$$
for fixed $\bb = (b_1,\ldots,b_n)$ and $F = \supp(\bb)$.  Summing over
all $\bb$ proves (\ref{eqn:1}).

The Cohen-Macaulay type of $S/M$ is $\sum_{\cc \in \NN^n} \beta_{d-1,
\cc} (M)$, so it suffices to prove that if $\beta_{d-1\,,\,\cc\,}(M) \neq
0$ then $\gamma_\cc \leq 1$.  Suppose the opposite, that is, $\gamma_\cc
\geq 2$ and $\beta_{d-1\,,\,\cc\,}(M) \neq 0$.  Then there are sets $F,F'
\subseteq {\{1, \ldots,n\}}$ such that $\cc \cdot F, \cc \cdot F' \in \irr
(S/M)$ are distinct.  Let $M_i = \mm^{\cc \cdot F}$ and $M_j = \mm^{\cc
\cdot F'}$ be the irreducible components $M$ corresponding to $\cc \cdot
F$ and $\cc \cdot F'$.  Since the algebraic co-Scarf complex of $M$ is
the minimal free resolution of $M$ and $\beta_{d-1\,,\,\cc\,}(M) \neq 0$,
there is an interior face $\sigma$ of the co-Scarf complex $\Delta_M^\aa$
with $\aa_\sigma = \DD - \cc$.  Since $m_i = \xx^{(\cc \cdot
F)^{\scriptstyle \aa}}$ and $m_j = \xx^{(\cc \cdot F')^{\scriptstyle
\aa}}$ divide $m_\sigma$ by construction, $\sigma$ contains both $i$ and
$j$.  In particular, $\{i, j\}$ is an edge of $\Delta_M^\aa$.  Now $S/M$ is
Cohen-Macaulay of codimension $\geq 2$, so $\supp (m_i) \cap \supp(m_j) \ne
\emptyset$ by Theorem~\ref{CMness2}.  But $\deg_{x_s} m_i = \deg_{x_s}
m_j = D - c_s >0$ for any $s \in \supp (m_i) \cap \supp(m_j)$,
contradicting the genericity of $M^\aa$. \qed
\vskip 2mm

After we had gotten the above proof, we conjectured the following more
general result about arbitrary triangulations of a simplex.  Margaret
Bayer proved our conjecture for {\it quasigeometric triangulations},
using local $h$-vectors \cite{Sta}. We are grateful for her permission
to include her proof in this paper.
Since the co-Scarf complex is a quasigeometric triangulation, 
Theorem \ref{marge} provides a second proof of Theorem~\ref{CMtype}.

\begin{thm}[M.~Bayer, personal communication] \label{marge}
Let $p_1$, $p_2$, \ldots, $p_r$ be points which lie in the relative 
interior of $(c-1)$-faces of a $(n-1)$-simplex $\Delta$.  Let $\Gamma$ be a 
quasigeometric triangulation of $\Delta$ having the $p_i$ among its 
vertices and having no interior $(n-c-1)$-face.  Then the number of 
interior $(n-c)$-faces is at least $r$.
\end{thm}

\begin{proof} 
According to the hypothesis, we have 
$\sum_{\substack{F \in \Delta \\ \#F=c}} f_0(\interior(\Gamma_F)) \ge r$, 
and  $f_i({\rm int}(\Gamma))=0$ for all $-1 \le i \le n-c-1$. 
By the decomposition of the $h$-polynomial of $\Gamma$ into local 
$h$-polynomials and the positivity of local $h$-vectors 
\cite[Theorem~4.6]{Sta}, we have 
$$h_{c-1}(\Gamma) \quad = 
\quad \sum_{F \in \Delta}\ell_{c-1}(\Gamma_F) \quad \ge \quad
\sum_{\substack{F \in \Delta \\ \#F=c}} \ell_{c-1}(\Gamma_F).$$ 
On the other hand, we have seen that  $\ell_{1}(\Gamma_F) = f_0(\interior
(\Gamma_F))$ in the proof of Theorem~\ref{size=c}. Since a local $h$-vector  
is symmetric \cite[Theorem~3.3]{Sta}, we have $\ell_{c-1}(\Gamma_F) = 
\ell_{1}(\Gamma_F) = f_0(\interior(\Gamma_F))$. So 
\begin{eqnarray*}
h_{c-1}(\Gamma)& \quad \ge & \quad
\sum_{\substack{F \in \Delta \\ \#F=c}} \ell_{c-1} 
(\Gamma_F) \quad = \quad
\sum_{\substack{F\in\Delta \\ \#F=c}} f_0(\interior(\Gamma_F))
\ge r. 
\end{eqnarray*}
Since the $h$-vector of $\interior(\Gamma)$ is the reverse of the $h$-vector of
$\Gamma$ (see the comment preceding \cite[Theorem~10.5]{Sta2}), we have 
\begin{eqnarray*}
h_{c-1}(\Gamma) &=& h_{n+1-c}(\interior(\Gamma))\\
&=& \sum_{i=0}^{n-c+1} (-1)^{n+1-c-i}\binom{n-i}{c-1} 
(f_{i-1}(\interior(\Gamma)))\\ 
&=& f_{n-c}(\interior(\Gamma)).
\end{eqnarray*}
Thus, the number of interior $(n-c)$-faces of $\Gamma$ is at least $r$.
\end{proof}

Our final results demonstrate the effective translation between
generic and cogeneric monomial ideals via Alexander duality.  

\begin{thm} \label{gens}
Let $M$ be a cogeneric monomial ideal with $r$ irreducible components,
each having the same codimension $c$. Then $M$ has at least 
$(c-1) \cdot r + 1$ minimal generators.
If $M$ has exactly $(c-1) \cdot r + 1$ generators then $S/M$ is 
Cohen-Macaulay. 
\end{thm}

\begin{proof}
The former statement is Alexander dual to Theorem~\ref{size=c}.  To prove
the latter statement, we recall the proof of Theorem~\ref{size=c}.
Assume that $S/M$ is not Cohen-Macaulay.  Then $\Gamma := \Delta_M^\aa$
has an edge $\{i,j\}$ whose excess $e$ satisfies $e \geq c$,
by Theorem \ref{CMness2}. Let $W \in
\Delta$ be the support of $m_{\{i,j\}}$.  Then $\#W = e+2$. By
\cite[Proposition~2.2]{Sta},
\begin{equation*} 
  \ell_W(\Gamma_W,x) =  \ell_2(\Gamma_W)x^2 + \ell_3(\Gamma_W) x^3 + \cdots,  
\end{equation*} where $\ell_2(\Gamma_W)$ is the number of edges of 
$\Gamma$ whose supports are $W$. So we have $f_{n-1}(\Gamma) = 
h(\Gamma,1) \geq (c-1) \cdot r + 1 + \ell_2(\Gamma_W) 
> (c-1) \cdot r + 1$ 
by an argument similar to the proof of Theorem~\ref{size=c}. 
Since $f_{n-1}(\Gamma)$ is equal to the number of generators of $M$, 
the proof is done. 
\end{proof}

\begin{exmp}\label{optimal}
The ideal  $M = \bigcap_{i=1}^r \< x_1^i, x_2^i, \cdots, x_{c-1}^i, 
x_{c-1+i}\>$ is cogeneric and has $(c-1) \cdot 
r+1$ minimal generators. Thus the inequality in Theorem~\ref{gens} 
is tight.
\end{exmp}

In the codimension $c=2$ case we can be more precise:

\begin{prop} \label{ht2}
Let $M$ be a cogeneric monomial ideal with $r$ irreducible components,
all of codimension $2$. Then $S/M$ is Cohen-Macaulay if and only if
$M$ has exactly $r + 1$ generators.
\end{prop}

\begin{proof} 
This is Alexander dual to Proposition~\ref{size=2}, in view of
Theorem~\ref{CMness2}.
\end{proof}

\vskip .5cm

\noindent{\bf Acknowledgements} \  We thank Margaret Bayer for contributing 
Theorem~\ref{marge} to this paper. Helpful comments on earlier drafts were  
given by Serkan Ho\c{s}ten and Irena Swanson. 
Work on this project was done during visits 
by Bernd Sturmfels to the Research Institute for Mathematical Sciences
of Kyoto University and by Kohji Yanagawa to the Mathematical Sciences 
Research Institute in Berkeley. The first two authors are 
also partially supported by the National Science Foundation.

\end{document}